\documentclass[12pt]{article}
\usepackage{amsfonts}
\begin{document}
\title{\bf  Variational problems related to some fractional kinetic
equations
\author{H. Hajaiej\\
 }}
\date{}
\maketitle
\begin{abstract}
We establish the existence and symmetry of all minimizers of a
constrained variational problem involving the fractional gradient.
This problem is closely connected to some fractional kinetic
equations.
\end{abstract}
\section{Introduction}

Fractional calculus has gained  a lot of interest during the last
decade due to its numerous applications in many fields. It appears
in wave propagation, inhomogenous porous material, geology,
hydrology, dynamics of earthquakes, bioengineering, chemical
engineering  signal processing, medicine, electrochemistry,
thermodynamics, neural networks, statistical physics, [2], [6], [8]
and references therein.\\
Fractional equations involving the fractional laplacian have also
played a crucial role in some kinetic problems, [5], [10], [11],
[12], [13] and [14], in which particular solutions are obtained by
solving the following minimization problem :
\begin{eqnarray*}
(P_c) : I_c &=& \inf \Big\{\int_{\mathbb{R}^N} |-\Delta^{s/2}(u)|^2
- \int_{\mathbb{R}^N}F(|x|, u) : u \in S_c\Big\}\\
E(u) &=& \frac{1}{2}\int_{\mathbb{R}^N} |-\Delta^{s/2}(u)|^2 -
\int_{\mathbb{R}^N}F(|x|,u),\\
F(r,t) &=& \int^t_0 f(x,p)dp,\\
S_c &=& \Big\{u \in H^s(\mathbb{R}^N) : \int_{\mathbb{R}^N} u^2 =
c^2\Big\},
\end{eqnarray*}
where $c$ is a prescribed number, $0 < s < 1$ and
$H^s(\mathbb{R}^N)$ is the usual Beso'v space , [8].\\
Note  also that the minimization problem $(P_c)$ appears in
disperative model equations : The generalized Benjamin-On equation,
the Benjamin-Bona-Mahong equation and the fractional nonlinear
Schr\"odinger equation.

In this paper, we   address the question of existence, radiality and
radial decreasiness of all minimizers of $(P_c)$ for integrands $F$
satisfying some growth conditions. This result generalizes a recent
one obtained by Frank and Lenzmann, [7], in which the authors have
considered the basic power nonlinearity $F(r,t) =
\frac{|t|^{\alpha+}}{\alpha+2}$.

Moreover, they have proved that the following minimization problems
$$I_M = \inf \left\{ \frac{1}{2} \int_{\mathbb{R}} |-\Delta^{s/2}u|^2 +
\frac{1}{\alpha+2} \int_{\mathbb{R}} |u|^{\alpha+2} : u \in
S_{\sqrt{M}}\right\} \eqno{(1.1)}$$ and
$$J^{s,\alpha} (Q) = \inf_{u \in H^s(\mathbb{R})\backslash \{0\}}
\frac{(\int_{\mathbb{R}}|-\Delta^{s/2}u|^2)^{\alpha/4s}
(\int_{\mathbb{R}}u^2)^{\frac{\alpha}{4s}(2s-1)+1}}{\int|
u|^{\alpha+2}}\eqno{(1.2)}$$
are equivalent.\\
They have also added that : $u$ is a solution of (1.1) if and only
if $u = e^{i\theta}\lambda^{1/\alpha} Q(\lambda^{1/2s}(.+y))$, for
some $\theta \in \mathbb{R}, y \in \mathbb{R}, \lambda > 0$ ; $Q$ is
a solution of (1.2).\\
Finally they have  stated that (1.2) (and therefore (1.1)) only has
minimizers when $0 < \alpha < \alpha_{\max}$, where $\alpha_{\max}$
is defined as follows :
$$\alpha_{\max} = \left\{ \begin{array}{ll}
\frac{4s}{1-2s}, 0 < s < \frac{1}{2}\\
\infty , \frac{1}{2} \leq s < 1
\end{array}\right.$$
This result seems to be erroneous and as we will show in section 2,
(1.1) admits minimizers if and only if $\alpha < 4s$.\\
Moreover, in this paper, we will study $(P_c)$ for general
nonlinearities $F$ such that : $|F(r,t)| \leq K (t^2 +
|t|^{\ell+2})$ where $0 < \ell < \frac{4s}{N}$.\\
Our main result, Theorem 2.1, states that :
\begin{enumerate}
\item If $0 < \ell < \frac{4s}{N}$, $(P_c)$ admits solutions and all
minimizers are radial and radially decreasing.
\item If $\ell = \frac{4s}{N}$, $(P_c)$ admits solutions and all
minimizers are radial and radially decreasing if $c^2$ is small
enough (some estimates will be given below).
\item If $\displaystyle{\liminf_{t\rightarrow \infty}}F(r,t)/t^{\ell+2} \geq A >
0$ for some $\ell > \frac{4s}{N}$, then $I_c = - \infty$ for all
$c$.\\
Now before stating our main result, let us first mention that
definitions and properties of the Schwarz symmetrization are
detailed in [4].
\end{enumerate}

If $u \in H^s_+(\mathbb{R}^N) = \{ u \in H^s(\mathbb{R}^N) : u \geq
0\}$, then the fractional Polya-Szeg\"o inequality holds true :
$$|\nabla_su^\ast|^2_2 = |-\Delta^{s/2}u^\ast|^2_2 \leq |-\Delta^{s/2}u|^2_2 =
|\nabla_su|^2_2\eqno{(1.3)}$$ which is a direct consequence of the
generalized Riesz inequality, [4], since as it was proven in [1] :
$$|\Delta^{s/2}u|^2_2 = C_{n,s} \int_{\mathbb{R}^N}\int_{\mathbb{R}^N}
\;\frac{|u(-x)-u(y)|^2}{|x-y|^{N+2s}}dx dy\;.$$ From now on $0 < s <
1, H^s(\mathbb{R}^N)$ is the standard Besov space. The norm of the
Lebesgue space $L^p(\mathbb{R}^N)$ is denoted by $|\,|_p$ . $c$ is a
prescribed number and $N \in \mathbb{N}^\ast$. In an integral where
no domain is given, it is to be understood that it extends on
$\mathbb{R}^N$.
\section{Main Result}
{\bf Theorem 2.1} Suppose that :\vspace{5mm}

$(F_0) F : [0,\infty) \times \mathbb{R} \rightarrow \mathbb{R}$ is a
Carath\'eodory function :

$\bullet\quad F(.,t) : [0, \infty) \rightarrow \mathbb{R}$ is
measurable in $\mathbb{R}_+ \backslash \Gamma$ for all $t \in
\mathbb{R}$, where $\Gamma$ is a subset of $\mathbb{R}_+$ having one
dimensional measure zero and ;

$\bullet\quad F(r,.) : \mathbb{R} \rightarrow \mathbb{R}$ is
continuous for every $r \in [0, \infty) \backslash \Gamma$.

$(F_1)\quad F(r,t) \leq F(r,|t|)$ for a.e $r \geq 0$ and every $t
\in \mathbb{R}$,\vspace{5mm}

$(F_2)$ For a.e $r \geq 0$ and every $t \geq 0$.
$$0 \leq F(r,t) \leq K(t^2 + t^{\ell+2}), \mbox{ where } K > 0 \mbox{ and }
0 < \ell < \frac{4s}{N}.$$

$(F_3)$ For every $\varepsilon > 0$, there exist $R_0 > 0$ and $t_0
> 0$ such that $F(r,t) \leq \varepsilon t^2$ for a.e $r \geq R_0$
and $0 \leq t \leq t_0$.

$(F_4) (v,y) \rightarrow F(\frac{1}{v} ,y)$ is supermodular on
$\mathbb{R}_+ \times \mathbb{R}_+$, i.e
$$F(r,a) + F(R,A) \geq F(r,A) + F(R,a)$$
for every $0 \leq r < R$ and $0 \leq a < A$.

Let $(\widetilde{P}_c) : \inf \{E(u) : u \in H^s_+(\mathbb{R}^N)$
and
$\displaystyle{\int_{\mathbb{R}^N}}u^2 \leq c^2\} = \widetilde{I}_c$\\
$(I) \quad\widetilde{I}_c < \widetilde{I}_d$ for $d < c$.\\
Then
\begin{enumerate}
\item $(P_c)$ admits a Schwarz symmetric minimizer for any $c$.
Moreover if $(F_4)$ holds with a strict sign, then for any $c$, all
minimizers of $(P_c)$ are radial and radially decreasing, i.e,
Schwarz symmetric.
\item If $\ell = \frac{4s}{N}$, then 1) holds true if and only if
$c$ is small enough.
\item If $\displaystyle{\liminf_{t\rightarrow \infty}}F(r,t)/t^\ell = A >
0$, with $\ell > \frac{4s}{N}$ then $I_c = - \infty$.
\end{enumerate}
{\bf Proof of 1} Fix $c$

{\bf Step 1} : {\bf$(P_c)$ is well  posed $(I_c > - \infty$ and all
minimizing sequences are bounded in} $H^s(\mathbb{R}^N))$.

By $(F_1)$ and $(F_2)$, we can write : $$\begin{array}{ll}
\int F(|x|, u(x))dx &\leq \displaystyle{\int} F(|x|, |u(x)|)dx\\
&\leq K c^2 + K \displaystyle{\int} |u(u(x)|^{\ell+2} dx
\end{array}\eqno{(2.1)}$$
Now using the fractional Gagliardo-Nirenberg inequality [8], [9], it
follows that : $$|u|_{\ell+2} \leq K'|u|^{1-\theta}_2|\nabla_s
u|^\theta_2\;\quad |\nabla_s u|_2 = \left(\int
|-\Delta^{s/2}u|^2\right)^{1/2} \mbox{ and } \theta =
\frac{N\ell}{2s(\ell+2)}.$$ Thus
$$\int |u(x)|^{\ell+2} dx \leq K'\{\int u^2(x)dx\}^{(1-\theta)(\ell+2)/2}
|\nabla_s u|^{\theta(\ell+2)}_2 . \eqno{(2.2)}$$
Therefore using Young inequality, we have\\
$ \{\displaystyle{\int}u^2(x)dx\}^{(1-\theta)
\frac{\ell+2}{2}}|\nabla_su|_2^{\theta(\ell+2)} \leq$
$$\frac{1}{p} \varepsilon^p \{|\nabla_s u|^2_2\}^{p\theta(\ell+2)/2}
+ \frac{1}{q\varepsilon^q} \left\{\int u^2
(x)dx\right\}^{q(1-\theta)\frac{(\ell+2)}{2}}\eqno{(2.3)}$$ for any
$\varepsilon > 0$ and $p > 1$ where $\frac{1}{p} + \frac{1}{q} = 1$.
\\
Choosing $p = \displaystyle{\frac{2}{\theta(\ell+2)} =
\frac{4s}{N\ell}}$, we get :\\
$\displaystyle{\int}|u(x)|^{\ell+2}dx \leq
\displaystyle{\frac{K'}{p}} \varepsilon^p\{|\nabla_s u|^2_2\} +
\frac{K'}{q\varepsilon^q}\{\displaystyle{\int}u^2(x)dx\}^{\frac{q(1-\theta)(\ell+2)}
{2}}$
$$= \frac{K'}{p} \varepsilon^p |\nabla_s u|^2_2 + \frac{K'}
{q\varepsilon^q} c^{q(1-\theta)(\ell+2)}
\eqno{(2.5)}$$ for any $u \in S_c$.\\
Therefore :  \begin{eqnarray*} E(u) &\geq& \frac{1}{2} |\nabla_s
u|^2_2 - K c^2 - K'K \varepsilon^p|\nabla_s u|^2_2-
\frac{KK'}{q\varepsilon^q} c^{q(1-\theta)(\ell+2)}\\
&=& (\frac{1}{2} - \frac{KK'}{p}\varepsilon^p)|\nabla_su|^2_2 -
Kc^2-\frac{KK'}{q\varepsilon^q} c^{q(1-\theta)(\ell+2)}.
\end{eqnarray*}
Thus $I_c > - \infty$ and all minimizing sequences are bounded  in
$H^s(\mathbb{R}^N)$.

{\bf Step 2} : {\bf Existence of a Schwarz symmetric minimizing
sequence}.\\
First note that if $u \in H^s(\mathbb{R}^N)$ then $|u| \in
H^s(\mathbb{R}^N)$.\\ Now by $(F_1)$, we certainly have :
$$E(|u|) \leq E(u)\quad \forall\; u \in H^s(\mathbb{R}^N).$$
Now by (1.3), we know that $|\nabla_s |u|^\ast|_2 \leq |\nabla_s
|u||_2$, and using Theorem 1 of [4], we have :
$$\int F(|x|, |u(x)|)dx \leq \int F(|x|, |u(x)|^\ast)dx$$
and
$$\int u^2 = \int(u^\ast)^2.$$
Thus without loss of generality, $(P_c)$ always admits a Schwarz
symmetric minimizing sequence.\vspace{5mm}\\
Let $(u_n)$ be a Schwarz symmetric minimizing sequence of $(P_c)$
for a fixed $c$.\vspace{5mm}

{\bf Step 3} : Let $(u_n) =(u^\ast_n)$ be a Schwarz symmetric
minimizing sequence then if $u_n $ converges weakly to $u$ $
\Rightarrow E(u) \leq \lim\inf
E(u_n)$.\\
{\bf Proof} $|\nabla_s u|_2 \leq \liminf |\nabla_s u_n|_2$ by the
weak lower semi-continuity of $\|\;\|_2$ of the fractional gradient
in
$H^s(\mathbb{R}^N)$.\\
Let us prove now that :
$$\lim_{n\rightarrow \infty} \int F(|x|, u_n(x)dx = \int F(|x|, u(x)).$$
Let $R > 0$, let us first prove that :
$$\lim_{n\rightarrow + \infty} \int_{|x|\leq R} F(|x|, u_n(x))dx =
\int_{|x| \leq R} F(|x|, u(x))dx.$$ Since $u_n$ converges weakly to
$u$ in $H^s(\mathbb{R}^N)$, it converges strongly  to $u$ in
$L^{\ell+2}(|x| \leq R)$. Thus there exists a subsequence
$(u_{n_k})$ of $(u_n)$ such that $u_{n_k} \rightarrow u$ a.e in
$L^2(B(0,R))$ and $|u_{n_k}| \leq h$ with $h \in L^{\ell+2}(|x| \leq
R). \{B(0,R) = \{x \in \mathbb{R}^N : |x| \leq R\}$\\
Now by $(F_2) : F(|x|, u_{n_k}(x)) \leq K(h^2(x)+h^{\ell+2}(x))$.\\
Noticing that $h^2 + h^{\ell+2} \in L^1(|x| \leq R)$, we get thanks
to the dominated convergence theorem :
$$\lim_{n\rightarrow + \infty} \int_{|x| \leq R}
F(|x|, u_n(x))dx = \int_{|x|\leq R}F(|x|, u(x))dx.$$ Let us prove
now that $\displaystyle{\lim_{R\rightarrow \infty}\lim_{n\rightarrow
\infty}\int_{|x|> R}}F(|x|, u(x))dx = 0$.\\
Let $n \in \mathbb{N}$, since $(u_n) = (u^\ast_n)$, we have that
$$V_N |x|^N u^2_n(x) \leq  \int_{|y| \leq |x|} u^2_n(y)dy \leq c^2.$$
Thus
$$u_n(x) \leq \frac{c}{V_N^{1/2}|x|^{N/2}} \leq \frac{c}{V_N^{1/2} R^{N/2}}\quad
\forall\; |x| > R.$$ Now let $\varepsilon  > 0$ and $R$ big enough,
we obtain thanks to $(F_3)$ that :
$$\int_{|x| > R}F(|x|, u_n(|x|) dx \leq \varepsilon \int_{|x|> R}
u^2_n(x)dx < \varepsilon c^2,$$ proving that
$$\lim_{R\rightarrow \infty}\lim_{n\rightarrow \infty}
\int_{|x| > R} F(|x|, u_n(x))dx = 0.$$ But $u$ inherits all the
properties of the sequence $(u_n)$ used to get the above limit, then
it follows that :
$$\lim_{R\rightarrow \infty}\int_{|x| > R} F(|x|, u(x))dx = 0.$$

{\bf Step 4} : {\bf $I_c$ is achieved}.\\
Denoting $v$ the weak limit of a Schwarz minimizing sequence of $(
\widetilde{P}_c)$. We certainly have, using previous steps, that
$$E(v) \leq \liminf E(u_n),$$
where
$$\lim_{n\rightarrow \infty} E(u_n) = \widetilde{I}_c. $$
On the other hand $|v|^2_2 = d^2 \leq c^2$.\\
It follows then by hypothesis $(I)$, that :
$$\widetilde{I}_c < \widetilde{I}_d \leq E(v) \leq \widetilde{I}_c$$
which is impossible, then $|v|^2_2 = c^2$. Suppose that $|v|^2_2 = d^2 < c^2$.\\
Therefore $I_c \leq E(v) = \widetilde{I}_c \leq I_c$, proving that
$(P_c)$ is achieved by $v = v^\ast$ a.e.

Now to show that all minimizers of $(P_c)$ are Schwarz symmetric, it
is sufficient to notice that if $(F_4)$ holds with a strict sign
then it follows by Theorem 1 of [4] that
$$E(u^\ast) < E(u)\quad \mbox{ for any } u \in H^s_t(\mathbb{R}^N)$$
and the result follows.\\
{\bf Proof of 2)} If $\ell = \displaystyle{\frac{4}{Ns}}$, (2.2)
becomes :
$$\int |u(x)|^{\ell+2} dx \leq K' c^{4/N}|\nabla_s u|^2_2\quad
\forall\; u \in S_c.$$ Hence
\begin{eqnarray*}
E(u) &\geq& \frac{1}{2} |\nabla_s u|^2_2 - K c^2 - KK' c^{4/N}
|\nabla_s u|^2_2\\
&=& (\frac{1}{2} -KK' c^{4/N})|\nabla_su|^2_2 - Kc^2.
\end{eqnarray*}
Thus $I_c > - \infty$ and all minimizing sequences are bounded in
$H^s(\mathbb{R}^N)$ provided that $0 < c <
(\frac{1}{2KK'})^{4/N}$.\\
Then previous steps (2,3 and 4) apply to $c$ such as in the latter
interval.\\
{\bf Proof of 3)} : It suffices to consider $u \in S_c$ and
$u_\lambda(x) = \lambda^{N/2}u(\lambda.) (\in S_c)$, then the
results follow when $\lambda$ tends to infinity.
\section*{References}
\begin{enumerate}
\item
  F. J . Almgren, E. H. Lieb : Symmetric decreasing is sometimes
continuous , J . Amer. Math. Soc ., 2 (1989) , 683-773.

\item G.Anastassiou : Advances on fractional inequalities, Springer
2011. \item  J. M. Angulo, V. V. Anh, R. Mc. Vinish, M .D.
Ruiz-Medina : Fractional kinetic equations driven by Gaussian  or
infinitely divisible noise. Advances in Applied Probability 2005.

\item A. Burchard, H Hajaiej: Rearrangement inequalities for
functional with monotone integrands. J of Functional Analysis, 233 (
2006), 561-582.

\item V B L Chaurasia, S C Pendy: Computable extensions of
generalized fractional kinetic equations in astrophysics. Research
in astronomy and astrophysics ( 2010).

\item  W Chen, S C Sun: Numerical solutions of fractional derivatives
equations in mechanics: advances and problems: The 22nd
international congress of theoretical and applied mechanics.
Australia 25-29 August 2008.

\item R. Frank, E. Lenzmann: Uniqueness  and non-degeneracy of ground
states  for $(-\Delta)^{s}Q+Q- Q^{\alpha+1}=0$ in $R$ , accepted in
Acta Mathematica.

\item H. Hajaiej, L Molinet, T Ozawa, B Wang : Necessary and
sufficient conditions for the fractional Gagliardo-Nirenberg
inequalities and applications to Navier-Stokes and generalized Boson
equations, Preprint.

\item H. Hajaiej, X. Yu, Z. Zhai : Fractional Gagliardo-Nirenberg and
Hardy inequalities  under Lorentz norms. Preprint.

\item T. A. M. Langlands, B. I.Henry : Fractional chemotaxis
diffusion equations. Physical review E- Statistical, Nonlinear and
soft matter physica (2010).

\item A. Mellet, S. Mischler, C. Mouhot : Fractional diffusion limit
for collisional kinetic equations, Archiv Rat Mechanical Analysis
(2008).

\item R. Nigmatullin : The 'Fractional' kinetic equations and general
theory of dielectric relaxation, vol 6, 5th international conference
on multibody systems, nonlinear dynamics and Control, Parts A, B, C
(2005).

\item R. K. Saxena, A. M. Mathai, H J Haubold : Solutions of  certain
fractional kinetic equations and a fractional diffusion equation.
Journal of Mathematical physics (2007).

\item R. K. Saxena, A. M. Mathai, H. J Haubold : On generalized
Fractional kinetic Equations. Physica A, Statistical Mechanics and
its applications ( 2004).
\end{enumerate}
\end{document}